\documentclass[a4paper, 8pt]{article}
\usepackage{graphicx}
\usepackage{amssymb}
\usepackage{latexsym, bm}
\usepackage{multicol}
\usepackage{indentfirst}
\usepackage{amssymb,amsfonts}
\usepackage{amsmath}
\usepackage{amsfonts}
\usepackage{fancyhdr}
\usepackage{titlesec}
\usepackage{cite}
\usepackage{ifthen}
\usepackage{pifont}
\usepackage{stmaryrd}
\usepackage{setspace}
\usepackage{indentfirst}
\usepackage{amsmath,amssymb,amscd,bbm,amsthm,mathrsfs,dsfont}
\usepackage{enumerate}
\newtheorem{theorem}{Theorem}[section]
\newtheorem{lemma}{Lemma}[section]

\newtheorem{definition}{Definition}[section]
\newtheorem{corollary}{Corollary}[section]
\newtheorem{claim}{Claim}[section]
\usepackage{amssymb}
\title{A note on the size Ramsey number of powers of paths
}
\author{Chunlin You\vspace*{3mm}
\\
{Center for Discrete Mathematics, Fuzhou University } \\
{Fujian 350108, China}\vspace*{3mm} \\
{\em (chunlin\_you@163.com)} }

\date{}
\begin{document}
\maketitle
\begin{abstract}
Let $r\geq3$ be an integer such that $r-2$ is a prime power and let $H$
be a connected graph on $n$ vertices with average degree at least $d$ and  $\alpha(H)\leq\beta n$, where $0<\beta<1$ is a constant. We prove that the size Ramsey number
\[
\widehat{R}({H};r) > \frac{{nd}}{2}{(r - 2)^2} - C\sqrt n
\]
for all sufficiently large $n$, where $C$ is a constant depending only on $r$ and $d$.
In particular, for integers $k\ge1$, and $r\ge3$ such that $r-2$ is a prime power, we have that there exists a constant $C$ depending only on $r$ and $d$ such that
$\widehat{R}(P_{n}^{k}; r)> kn{(r - 2)^2}-C\sqrt n -\frac{{({k^2} + k)}}{2}{(r - 2)^2}$ for all sufficiently large $n$, where $P_{n}^{k}$ is the $kth$ power of $P_n$.
We also prove that $\widehat{R}(P_n,P_n,P_n)<764.1n$ for sufficiently large $n$. This result improves some results of Dudek and Pra{\l}at (\emph{SIAM J. Discrete Math.}, 31 (2017), 2079--2092 and \emph{Electron. J. Combin.}, 25 (2018), no.3, \# P3.35).


\end{abstract}

\section{Introduction}
For a graph $G$, we use $V(G)$ and $E(G)$ to denote its vertex set and edge set, respectively.
Denote $|E(G)|=e(G)$ and $|V(G)|=v(G)$.
For a vertex $v$ in $V(G)$, the degree of $v$ is denoted by $d_G(v)$,
which is the number of vertices of $G$ that adjacent to $v$.
We write $K_N$ for a complete graph on $N$ vertices.
For a subset $I\subseteq V(G)$, we say $I$ is an independent set if any two vertices of $I$ are not adjacent.
The independence number of $G$, denoted $\alpha(G)$, is the order of the largest
independent set in $G$.

Given a graph graphs $G$, $H$  and a positive integer $r$,
the Ramsey number $R(H;r)$ of $H$ is the smallest
number $N$ such that any $r$-coloring of $E(K_N)$ contains a monochromatic copy of $H$.
Instead of minimizing the number of vertices, one can minimum number of edges.
This naturally leads to the size Ramsey number $\widehat{R}({H;r})$ introduced by Erd\H{o}s, Faudree, Rousseau and Schelp \cite{Erd-R-1978}.
The size Ramsey number $\widehat{R}({H;r})$ is the minimum integer $m$ such that there exists a graph
$G$ on $m$ edges such that every coloring of the edges of $G$ with $r$ colors yields
a monochromatic copy of $H$.
When $r=2$, we denote $\widehat{R}({H})$ instead by $\widehat{R}({H;r})$ for short.

Note that we have the trivial upper bound
$\widehat{R}({H;r}) \leq {{R}({H;r})\choose 2}$ for any graph $H$.
When $H$ is complete graph, Erd\H{o}s, Faudree, Rousseau and Schelp \cite{Erd-R-1978}
proved that
$\widehat{R}({K_n}) = {R({K_n})\choose 2}$.
In \cite{Erd-c},  Erd\H{o}s proved that
$\frac{1}{60}n^{2}2^{n}<\widehat{R}(K_{n,n})<\frac{3}{2}n^{3}2^{n}.$
For the two color size  Ramsey number of the path $P_{n}$ on $n$ vertices, Erd\H{o}s \cite{Erd-1981}
offered 100 dollars for a proof or disproof that
$\widehat{R}({P_n}))/n \to \infty  ~~ \text{and} ~~ \widehat{R}({P_n})/{n^2} \to 0.$
In 1983, Beck gained the 100 dollars prize by proving that $\widehat{R}({P_n}) < 900n $ for sufficiently large $n$.
Alon and Chung \cite{Alon-1988} provided an explicit construction of a graph $G$ with $O(n)$ edges such that $G\to P_n$.
After many successive improvements (\cite{boll-1986,dudek-2017}) for the lower bound,
\cite{Boll-2001, dudek-2015, letzer-2016, dudek-2017, bal-2019} (for the upper bound) the state of the art is
$(3.75-o(1))n \leq \widehat{R}({P_n})\leq 74n.$
For any $r\geq 2$ colors, Dudek and Pra{\l}at \cite{dudek-Pr-2018}
proved
$\widehat{R}({P_n};r)\leq Cr^{2}(\log r)n$.
As for the lower bound,  Dudek and Pra{\l}at \cite{dudek-2017} proved that for any $r\geq 2$,
$\widehat{R}({P_n};r)\geq  \frac{{(r + 1)(r+2)}}{4}n$ and then
Krivelevich\cite{Pokrovskiy-2018} proved that for $r\geq 3$ and $r-2$ is a prime power, $\widehat{R}(P_n;r)>(r-2)^{2}n-O(\sqrt{n})$ for all sufficiently large $n$.
Recently, Bal and DeBiasio \cite{bal-2019} improve on each of these results by proving the following:
\[\widehat{R}({P_n};r) \ge \max \left\{ {\left( {\frac{{(r - 1)r}}{2} + 2.75 - o(1)} \right)n,({q^2} - o(1))n} \right\}\]
for $r\geq 2$, where $q$ is the largest prime power such that $q\leq r-1$.

Moving away from paths, Beck \cite{beck-1990} even asked if there is some constant $c=c(\Delta)>0$ such that
$\widehat{R}(G)\leq cn$ for any graph $G$ with $n$ vertices and
maximum degree at most $\Delta$.
R\"{o}dl and Szemer\'{e}di \cite{V.Rodl-Eo(2000)} proved
that there exists positive constants $c$ and $\alpha$, and a graph $G=(V,E)$
with $|V|=n$ and maximum degree $\triangle(G)=3$ such that
$\widehat{R}(G)\geq cn(\log_{2}n)^{\alpha}$.
The current best known upper bound for the class of all graphs with constant maximum degree
was due to Kohayakawa, R\"{o}dl, Schacht and Szemer\'{e}di \cite{Kohayakawa-rodl-schacht-2011}.
For every graph $G$ on $n$ vertices of maximum degree $d$,
\[
\widehat{R}(G) \le c{n^{2 - 1/d}}{\log ^{1/d}}n.
\]
For more results on the size Ramsey number of bounded degree graphs see
\cite{Domingos, Haxell-Kohayakawa-Luzak, Kohayakawa-Retter-2019, Clemens-2019, Pokrovskiy-2017, sun-li}.
See also \cite{Ben-Eliezer-2012, Balogh-Clemen-2018} for related results in a digraph setting.

Let us turn our attention to powers of bounded degree graphs.
Let $H$ be a graph with $n$ vertices and $k\geq 1$ be an integer, the $kth$ power of $H$, denote
$H^{k}$, is the graph with vertex set $V(H)$ where distinct vertices $u$ and
$v$ are adjacent iff the distance between them in  $H$ is at most $k$.
Recently it was proved that the 2-color size Ramsey number
of powers of paths and cycles is linear \cite{Clemens-Jenssen-Kohayakawa}.
This result was extended to any fixed number $r$ of
colours in \cite{han-jenssen-2018}, i.e.,
\begin{equation}\label{inequ-1}
\begin{split}
\widehat{R}(P_{n}^{k};r) = O(n).
\end{split}
\end{equation}
Recently,   (\ref{inequ-1}) was extended  to  powers of bounded degree trees in \cite{Berger-2019}.

In this paper,  we prove a nontrivial lower bound of the multicolor size Ramsey numbers for graph with bounded average degree and independence number.
As a corollary, we give a nontrivial lower bound on (\ref{inequ-1}).

\begin{theorem}\label{restrict degree}
Let $r\geq3$ be integer such that $r-2$ is a prime power, and let $H$
be a connected graph on $n$ vertices with average degree at least $d$ and  $\alpha(H)\leq\beta n$,  where $0<\beta<1$ is a constant. We have
\[
\widehat{R}({H};r) > \frac{{nd}}{2}{(r - 2)^2} - C\sqrt n
\]
for all sufficiently large $n$, where $C=C(r,d)$ is a constant depending only on $r$ and $d$.
\end{theorem}

Recently, Han et al. \cite{han-jenssen-2018} proved that the multicolor size Ramsey number of powers of paths $P_{n}^{k}$ is linear.
Since their proof is based on the regularity lemma and so no specific constant is known.
Based on our main result above, we can give a specific lower bound for $\widehat{R}(P_{n}^{k}; r)$.
Note that the average degree of $P_{n}^{k}$ is $d(P_{n}^{k})=\frac{{2nk - {k^2} - k}}{n}$.
The following corollary is immediate from Theorem \ref{restrict degree}.

\begin{corollary}
Let $r\geq3$ be integer such that $r-2$ is a prime power.
\[
\widehat{R}(P_{n}^{k}; r)> kn{(r - 2)^2} -\frac{{({k^2} + k)}}{2}{(r - 2)^2}-C\sqrt n
\]
for all sufficiently large $n$ and $k\geq 2$, where $C=C(r,d)$ is a constant depending only on $r$ and $d$.
\end{corollary}

Next,  we  give an upper bound for $\widehat{R}(P_n,P_n,P_n)$ as follows, which is better than that by Dudek and Pra{\l}at \cite{dudek-2017,dudek-Pr-2018} and that by Krivelevich \cite{Pokrovskiy-2018}.
\begin{theorem}\label{th-p3}
Let $c=8.2919$ and $d=82.1405$. Then, a.a.s ${g_{cn,cn,d}} \to (P_n)_{3}$, which implies that
$\widehat{R}(P_n,P_n,P_n)<764.1n$ for sufficiently large $n$.
\end{theorem}

\section{Proof of Theorem \ref{restrict degree}} \label{chap2.3}

Our argument extends the methods of Krivelevich \cite{Pokrovskiy-2018}.
We first introduce the affine plane.
In geometry, an affine plane is a system of points and lines that satisfy the following properties (Let $A_{q}$ denotes an affine plane of order $q$):
\begin{enumerate}[(1)]
\item There are $q^{2}$ points.
  \item Every line has $q$ points.
  \item Any pair of distinct points is contained in a unique line.
  \item There are $q^2+q$ edges (lines).
  \item These lines can be split into $q+1$ disjoint families (ideal points), such that the lines in every ideal point do not intersect.
\end{enumerate}
More information refer to  Krivelevich \cite{Pokrovskiy-2018}. Such a system is known
to exist for every prime power $q$.
We give an example when the affine plane of order $q=3$.
see the following Fig \ref{fig:Affine_plane}. We give these points label figures.

\begin{figure}
\begin{center}
 \includegraphics[scale=0.4]{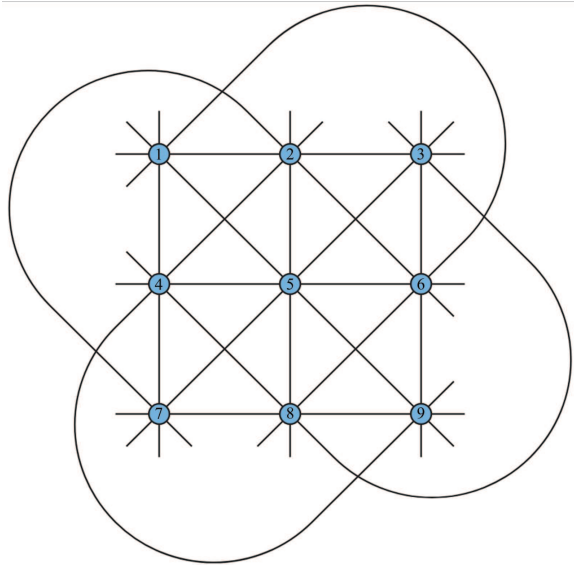}
\end{center}
  \caption{\small An affine plane of order $q=3$.}
  \label{fig:Affine_plane}
\end{figure}

When $q=3$ the ideal points including four sections:
\begin{itemize}
  \item {\textcircled{1}\textcircled{2}\textcircled{3}}; {\textcircled{4}\textcircled{5}\textcircled{6}}; {\textcircled{7}\textcircled{8}\textcircled{9}}.
  \item \textcircled{1}\textcircled{4}\textcircled{7}; \textcircled{2}\textcircled{5}\textcircled{8}; \textcircled{3}\textcircled{6}\textcircled{9}.
  \item \textcircled{1}\textcircled{5}\textcircled{9}; \textcircled{2}\textcircled{6}\textcircled{7}; \textcircled{3}\textcircled{4}\textcircled{8}.
  \item \textcircled{3}\textcircled{5}\textcircled{7}; \textcircled{1}\textcircled{6}\textcircled{8}; \textcircled{2}\textcircled{4}\textcircled{9}.
\end{itemize}
At each section different line do not intersect, these section is ideal points.

We will use the bounded difference inequality (BDI).
This is a concentration inequality that generalizes Hoeffding's inequality and
has found many uses in learning theory.
The BDI was first proved by McDiarmid  \cite{McDiarmid}.

\begin{definition}{\emph{(Bounded difference assumption)}}.\label{BDI-DEF}
Let $A$ be some set and
\[
\varphi :{A^n} \to R.
\]
We say $\varphi$ satisfies the bounded
difference assumption if $\exists\; {a_1},{a_2},\dots,{a_n} \ge 0 ,$
such that for any  $1 \le i \le n$
\[
{\mathop {\sup }\limits_{x_1,x_2, \ldots ,x_n,x_{i}' \in A} \left| {\varphi (x_1, \ldots ,x_i, \ldots ,x_n) - \varphi (x_1, \ldots ,x_i, \ldots ,x_n)} \right| \le a_i}
.\]
That is, if we subsitute  $x_i$ to $x_{i}'$, while keeping other $x_j$
fixed. $\varphi$ changes by at most $a_i$.
\end{definition}

\begin{theorem}[McDiarmid \cite{McDiarmid}]\label{BDI}
Let $X_1,X_2,\dots,X_n$ be arbitrary independent random variables on set $A$
and $\varphi:{A^n} \to R$ satisfy the bounded difference assumption.
Then $\forall \gamma  > 0$,
we have
\[
\Pr\left\{\varphi(X_1,X_2,\dots,X_n)-E[\varphi(X_1,X_2,\dots,X_n)]\geq \gamma\right\}
\leq {\exp\left\{ - {{2{\gamma ^2}}}/{{\sum\limits_{i = 1}^n {{a_i}^2} }}\right\}}
.\]
\end{theorem}

\bigskip
Now, we are ready to give a proof of our main result.

\bigskip
\noindent{\bf Proof of Theorem \ref{restrict degree}.}
Let $G=(V,E)$ be a graph with
\[
\left| E \right| \le \frac{{nd}}{2}{(r - 2)^2} - C\sqrt n
\]
edges. Without loss of generality, we can assume that $G$ is connected.

We aim to prove that there exists an $r$-coloring $C_1$ of $E(G)$ such that it contains no monochromatic $H$ whose average degree is at least $d$ and  $\alpha(H)\leq \beta n$,  where $0<\beta<1$ is a constant.

Given $q=r-2$, and let $A_{q}$ be an affine plane of order $q$, where $q$ is a prime power.
We assume that the vertex  set of $A_q$ is $\{1,2,\dots,q^2\}$, and the ideal points
are $I_1,I_2,\dots,I_{q+1}$.

We now split vertex set $V$ of $G$ into two sets $V_{0}$ and $V\setminus V_{0}$, where
\[
V_{0}=\left\{v\in V(G): d_{G}(v)\geq \frac{1}{{1 - \beta }}r^{2}d\right\}.
\]
Clearly, $|V_{0}|<(1-\beta)n$.
Partition $V\setminus V_0$ at random into $q^{2} $ parts $V_{1},V_{2},\dots,V_{q^{2}}$
by puting every vertex $v \in V\backslash {V_0}$ into $V_{i}$, $1 \le i \le {q^2}$,
independently and with probability $1/q^2$.
Fix a partition $(V_{1},V_2,\dots,{V_{{q^2}}})$, we define a coloring $C_1$ as follows:
\begin{enumerate}[(1)]
  \item If the endpoints of $e$ fall into $V_x$ and $V_y$ for $1\leq x <y\leq q^2$ , $L$ is the unique line of $A_q$ passing through $x$,
  $y$ then $C_{1}(e)=i$, whenever $L\in I_i$.

  \item If $e$ falls inside one of the $V_x$, then color it an arbitrary color different from $r$.

  \item If at least one of end vertices of $e$ falls inside $V_{0}$, then $C_{1}(e)=r.$
\end{enumerate}

For convenience, denote $G_i$ by the subgraph induced by all edges in color $i$, where $1\le i\le r$.

It is clear that $G_r$ does not contain $H$. If not, since the edges in color $r$ are only those incident to $V_{0}$, it follows that the subgraph induced by edges in color $r$ has independence number at least $|V\setminus V_0|>\beta n$. But this contracts the $\alpha(H)\leq\beta n$.

For color $i$, $1\leq i \leq q+1=r-1$,
we have the following claim.

\begin{claim}
Each connected component of $G_i$ can only contain those edges corresponding to vertex subset $\bigcup\limits_{x \in L} {{V_x}} $
for some line $L\in I_i$.
\end{claim}
\noindent
{\bf Proof.} Fixed a connected component of $G_i$, which we denote it by $D_i$.
Suppose to the contrary, there exists an edge $e=uv$ in $D_i$
such that the unique line $L_1$ passing through $V_1$ and $V_2$ with $u\in V_1$ and $v\in V_2 $ which is distinct from $L$.
According to the  above  coloring  method,  we have $L_1 \in I_i$.
Also, since the lines of the ideal point are pairwise disjoint, which means
\[
L \cap {L_1} = \emptyset.
\]
Thus  $V_1$ and $V_2$ are distinct from $V_x$ for ${x \in L}$, and
\[
u,v \notin \bigcup\limits_{x \in L} {{V_x}}.
\]

Note that $D_i$ is connected, there must exist $w\in V_3\subseteq\cup_{x \in L} {{V_x}}$ satisfying
$uw\in E(D_i)$ or $vw\in E(D_i)$. Without loss of generality, assume that $uw\in E(D_i)$.
Suppose $L_2$ is the unique line of $A_q$ passing through $V_{1}$ and $V_{3}$.
However, $C_1(u,w)=i$ and hence $L_2\in I_i$, this clearly leads to a contradiction that $u\in L_1\cap {L_2}\ne \emptyset$ from the fact that the lines of the ideal point are pairwise disjoint.
In conclusion, the proof of the claim is now complete.
\hfill $\Box$

Now, it remains to check that each line $L$ of $A_q$ corresponding to vertex subset $\cup_{x \in L} {{V_x}} $ spans fewer than $nd/2$ edges.

Consider an edge $e\in E(G) $. For a given line $L$ of $A_q$,
the probability that $e\subseteq \cup_{x \in L} {{V_x}}$ is
\[
{\left(\frac{{| L |}}{{{q^2}}}\right)^2}=\frac{1}{{{q^2}}} = \frac{1}{{{{(r - 2)}^2}}}.
\]

Let $A_L$ be the random variable such that
${A_L} = | {E(\bigcup\limits_{x \in L} {{V_x}} )} |$.
According to the linearity of expectation,
it follows that
\[E[{A_L}] = | {E(G)} |.\frac{1}{{{{(r - 2)}^2}}} \le \frac{{nd}}{2} - \frac{C}{{{{(r - 2)}^2}}}\sqrt n <nd/2. \]

In order to complete the proof of Theorem \ref{restrict degree}, it suffices to prove the random variables $A_L$ are concentrated.

According to the previous proof process.
Observe that changing the location (part $V_x$) of a vertex $v\in V\setminus{V_0}$
changes the count in any $A_L$ by at most $\frac{1}{{1 - \beta }}r^{2}d$.
From Definition \ref{BDI-DEF}, we know the $a_i=\frac{1}{{1 - \beta }}r^{2}d$.
In order to use Theorem \ref{BDI}, we set
$\varphi(X_1,X_2,\dots,X_n)=A_L$ and $\gamma=\frac{C}{{{{(r - 2)}^2}}}\sqrt n$.

Thus we have
\[
\Pr\left\{A_L-E[A_L]\geq \gamma\right\}\leq {\exp\left\{\frac{{ - 2{C^2}{{(1 - \beta )}^2}}}{{{r^8}{d^2}\beta }}\right\}}
.\]

Applying the union bound over all lines $L$ of $A_q$ and taking $C=C(r,d)$
to be large enough,
we conclude that
\[
\Pr \left\{ {\mathop  \bigcup \limits_L \left({A_L} - E\left[ {{A_L}} \right] < \gamma \right)} \right\} > 0,
\]
where $\gamma=\frac{C}{{{{(r - 2)}^2}}}\sqrt n$.
That is ${A_L}<E[{A_L}]+\gamma<nd/2$ have positive probability for any $L$.

Thus, the random variables $A_L$ are concentrated.
In conclusion, this completes the proof of Theorem \ref{restrict degree}. \hfill$\Box$

\section{Proof of Theorem \ref{th-p3}}\label{chap4}
We will focus on upper bounds for size Ramsey number of paths more than two colors.
Here is a natural generalization of pairing  model in random bipartite graph.
We use the \textit{pairing model}  on the  random regular graphs.
It was first introduced by Bollob\'{a}s \cite{boll-1980}.
In this section, it is extended to the bipartite graph.

Now let's briefly introduce this method.
We divide $2cdn$ points into $2cn$ boxes. These boxes are labeled $v_1,v_2,\dots,v_{2cn}$,
then each box has $d$ points.
A \textit{pairing}
of these points is a perfect matching into $cdn$ pairs. Given a pairing $P$, we can construct a
$d$-regular multigraph $G(P)$.
The vertex set of $G(P)$ is the boxes: $v_1,v_2,\dots,v_{2cn}$.
If  $v_iv_j\in E(G(P))$ owing to a pair $\{ {x,y}\}$ among $x\in v_i $  and  $y\in v_j$.

For a model of random bipartite regular,
we consider graphs with $2cn$ vertices and $V=\{1,2,\dots,2cn\}$ to be the vertex
set.
we assume that the vertices of one color
are labelled $1,2,\dots,cn$ as are the vertices of the other color.
The natural pairing model for $d$-regular bipartite graphs is obvious.
The boxes containing points
are the same as for ordinary graphs.
However, the random perfect matching is equivalent to a bijection
between the points in boxes of one colour and those in boxes of the other
colour.
So, the bipartite graph is easier than in the graph case.

In this case, the loops and parallel edges allowed.
We need to limit this random bipartite graph to a simple bipartite graph.
O'Neil \cite{neil} found that the number of bicoloured $d$-regular graph on $n$
vertices ($n$ even) is asymptotic to
\[
\frac{{(\frac{1}{2}dn)!{e^{ - \frac{1}{2}{{(d - 1)}^2}}}}}{{{{(d!)}^n}}},
\]
for $3 \le d \le {(\log n)^{\frac{1}{4} - \varepsilon }}$ with $\varepsilon>0$.
We know that the probability of a simple bipartite regular graph by this random  pairing is
\[
{{e^{ - \frac{1}{2}{{(d - 1)}^2}}}}.
\]
This probability only depends on $d$.
So the graph generated by this probability space is a.a.s a simple graph.
For more information on this model,
see, for instance, the survey of Wormald \cite{wormald-1999}.

Before giving a proof for Theorem \ref{th-p3}, we need the following  useful lemma by  Dudek and Pra{\l}at \cite[Lemma 3.7]{dudek-2017}.

\begin{lemma}[Dudek and Pra{\l}at \cite{dudek-2017}]\label{lem-dudek}
Let $r \ge2$ and $G =( {{V_1} \cup {V_2},E})$
be a balanced  bipartite graph of order $cn$ for some $c>2^{r}-1$. Assume that for any subsets $S\subseteq V_1$ and $T\subseteq V_2$ with $|S|=|T|=((c+1)/2^{r}-1)n/2$, we have $e(S,T)\neq 0$. Then, $G\rightarrow (P_n)_r$.
\end{lemma}

We will use the following form of Lemma \ref{lem-dudek}.
\begin{lemma}\label{p-3}
Let $G =( {{V_1} \cup {V_2},E})$ be a balanced bipartite graph on $2cn$ ($c>3.5$) vertices. Assume that for
any subsets $S\subseteq V_1$ and $T\subseteq V_2$ with $|S|=|T|=n(2c - 7)/16$, we have $e(S,T)\neq 0$. Then, $G\rightarrow (P_n)_{3}$.
\end{lemma}

\noindent
{\bf Proof of Theorem \ref{th-p3}}.
Consider $g_{cn,cn,d}$ for some $c\in(3.5,\infty )$ and $d\in N$. Our goal is to show that, for some suitable
choice of $c$ and $d$, the expected number of pairs of two disjoint sets $S$
and $T$, where $ S \subseteq {V_1}$ and $ T \subseteq {V_2}$ such that $|S|=|T|=n(2c-7)/16$ and  $e(S,T)=0$
tends to zero as $n\rightarrow \infty $.
This together with the first moment principle, implies that a.a.s. no such pair exists and so,
by Lemma \ref{p-3} , we get that a.a.s. ${g_{cn,cn,d}} \to (P_n)_{3}$.  As a result,
$\widehat{R}(P_n,P_n,P_n)\leq(cd+o(1))n $.

Let $c_{1}=(2c-7)/16$ and $X(c,d)$ be the expected number of pairs of two disjoint sets $S,T$ such that $|S|=|T|=n(2c-7)/16$,  $e(S,T)=0$, and $ e(S,{V_2}\backslash T) = c_{1}dn $, $ e(T,{V_1}\backslash S) = c_{1}dn $.
Using the pair model, it is clear that we have

\begin{align*}
\nonumber X(c,d)=&\binom{cn}{c_{1}n}\binom{cn}{c_{2}n}\binom{(c-c_{1})dn}{c_{1}dn}\binom{(c-c_{1})dn}{c_{1}dn}(c_{1}dn)!(c_{1}dn)!
\\& \times M\big((c-2c_{1})dn, (c-2c_{1})dn \big)/M(cdn, cdn),
\end{align*}
where $M(t,t)$ is the number of perfect matchings on balanced  bipartite graph of order $2t$ vertices, that is,
$M(t,t)=t!.$

Using Stirling's formula ($t! \sim \sqrt {2\pi t} {(t/e)^t}$), after simplification we get

\[
X(c,d) = f(c,d){e^{g(c,d)n}},
\]
where
\[
f(c,d) = \frac{1}{{2\pi }}\cdot\frac{1}{{{c_1}n}}\sqrt {\frac{c}{{c - 2{c_1}}}}
\]
and
\begin{align*}
g(c,d)&=2c\ln c + 2(c-c_{1})d\ln (c-{c_1})-2{c_1}\ln{c_1}
\\&- 2(c-c_1)\ln(c-c_1)
-(c-2c_{1})d\ln(c-2c_{1})-cd\ln c.
\end{align*}

It remains to choose suitable $c$ and $d$ such that $g(c,d)\leq 0$ and $cd$ as small as possible.
By using Matlab, we can take $c=8.2919$ and $d=82.1405$, it follows that
\[
\widehat{R}(P_n,P_n,P_n)\le cdn<764.1n,
\] for sufficiently large $n$. This completes the proof of Theorem \ref{th-p3}.
\hfill $\Box$

\subsection{Concluding Remark}\label{remark}
Indeed, we can use the same method as in Theorem \ref{th-p3} to get the upper bounds for $r=4,5$
that are listed in Table 1 as follows, in which  $n$ is sufficiently large.
We omit the proofs of these results since they are quite similar to that of Theorem \ref{th-p3}.

\begin{small}
\begin{center}
 \setlength{\tabcolsep}{8mm}
{\bf Table 1.} Comparing upper bounds with Dudek and Pra{\l}at \cite{dudek-2017,dudek-Pr-2018}

\medskip
\begin{tabular}{|c|c|c|c|c|c|c|}  \hline
$r$  & $3$ & $4$ & $5$ \\ \hline
        $\widehat{R}(P_n;r)/n<$     & 764.1  & 5167.7   &  56110    \\   \hline
         \cite{dudek-2017,dudek-Pr-2018} $\widehat{R}(P_n;r)/n<$     & 6336, 5933  & 33792, 13309  &168960, 24142    \\   \hline
       \end{tabular}
\end{center}
\end{small}

From  Table 1, we can see that our results are better than  Dudek and Pra{\l}at for
$r=3, 4$.

\end{document}